\pdfoutput=1 
\documentclass[11pt]{amsart}

\usepackage{epsfig,overpic,comment}
\usepackage[usenames,dvipsnames,svgnames,table]{xcolor}
\usepackage[hyphens]{url}
\usepackage[pagebackref,linktocpage=true,colorlinks=true,linkcolor=Blue,citecolor=BrickRed,urlcolor=RoyalBlue]{hyperref}
\usepackage[msc-links,abbrev]{amsrefs}
\usepackage{amsmath,amsthm,amssymb,marginnote}
\usepackage{enumerate}
\usepackage[T1]{fontenc}

\newtheorem{theorem}{Theorem}[section]

\theoremstyle{definition}

\newcommand{\be}{\begin{equation}}
\newcommand{\ee}{\end{equation}}
\newcommand{\ol}{\overline}

\newcommand{\R}{\mathbf{R}}

\newcommand{\C}{\mathcal{C}}
\newcommand{\G}{\Gamma}
\newcommand{\M}{\mathcal{M}}

\renewcommand{\epsilon}{\varepsilon}

\makeatletter
\DeclareFontFamily{U}{tipa}{}
\DeclareFontShape{U}{tipa}{m}{n}{<->tipa10}{}
\newcommand{\arc@char}{{\usefont{U}{tipa}{m}{n}\symbol{62}}}%

\newcommand{\arc}[1]{\mathpalette\arc@arc{#1}}

\newcommand{\arc@arc}[2]{%
  \sbox0{$\m@th#1#2$}%
  \vbox{
    \hbox{\resizebox{\wd0}{\height}{\arc@char}}
    \nointerlineskip
    \box0
  }%
}
\makeatother

\makeatletter
\def\@tocline#1#2#3#4#5#6#7{\relax
  \ifnum #1>\c@tocdepth 
  \else
    \par \addpenalty\@secpenalty\addvspace{#2}%
    \begingroup \hyphenpenalty\@M
    \@ifempty{#4}{%
      \@tempdima\csname r@tocindent\number#1\endcsname\relax
    }{%
      \@tempdima#4\relax
    }%
    \parindent\z@ \leftskip#3\relax \advance\leftskip\@tempdima\relax
    \rightskip\@pnumwidth plus4em \parfillskip-\@pnumwidth
    #5\leavevmode\hskip-\@tempdima
      \ifcase #1
       \or\or \hskip 1.3em \or \hskip 2em \else \hskip 5em \fi%
      #6\nobreak\relax
    \hfill\hbox to\@pnumwidth{\@tocpagenum{#7}}\par
    \nobreak
    \endgroup
  \fi}
\makeatother

\newcommand{\nocontentsline}[3]{}
\newcommand{\tocless}[2]{\bgroup\let\addcontentsline=\nocontentsline#1{#2}\egroup}

\begin{document}
\setlength{\baselineskip}{1.2\baselineskip}

\title[Total mean curvatures of Riemannian hypersurfaces] 
{Comparison formulas for total mean curvatures\\ of Riemannian hypersurfaces}

\author{Mohammad Ghomi}
\address{School of Mathematics, Georgia Institute of Technology,
Atlanta, GA 30332}
\email{ghomi@math.gatech.edu}
\urladdr{www.math.gatech.edu/~ghomi}

\vspace*{-0.75in}
\begin{abstract}
We devise some differential forms after Chern to compute a family of formulas 
for comparing total mean curvatures of nested hypersurfaces in Riemannian manifolds. This yields a quicker proof of a recent result of the author with Joel Spruck, which had been obtained via Reilly's identities. 
\end{abstract}

\date{\today \,(Last Typeset)}
\thanks{The research of the author was supported by NSF grant DMS-2202337.}

\maketitle


\section{Introduction}\label{sec:intro}
The  \emph{total $r^{th}$  mean curvature} of an oriented $\C^{1,1}$  hypersurface $\G$ in a Riemannian $n$-manifold $M$, for $0\leq r\leq n-1$,  is given by
$$
\M_{r}(\G):=\int_\G \sigma_r(\kappa), 
$$
where 
$\kappa:=(\kappa_1,\dots,\kappa_{n-1})$ denotes  the principal curvatures of $\G$, with respect to the choice of orientation, and $\sigma_r\colon \R^{n-1}\to\R$ is the $r^{th}$ symmetric function; so
$$
\sigma_r(\kappa)=\sum_{1\leq i_1<\dots<i_r\leq n-1}\kappa_{i_1}\dots \kappa_{i_r}.
$$
We  set $\sigma_0:=1$, and $\sigma_r:=0$ for $r\geq n$ by convention. Thus $\mathcal{M}_0(\Gamma)$ is the $(n-1)$-dimensional volume, $\mathcal{M}_1(\Gamma)$ is the total mean curvature, and $\mathcal{M}_{n-1}(\Gamma)$ is the total Gauss-Kronecker curvature of $\Gamma$. Up to  multiplicative constants, these quantities form the coefficients of Steiner's polynomial, and are known as quermassintegrals when $\Gamma$ is a convex hypersurface in Euclidean space.
The following result was established in \cite[Thm. 3.1]{ghomi-spruck2023} generalizing earlier work in \cite[Thm. 4.7]{ghomi-spruck2022}:

\begin{theorem}[\cite{ghomi-spruck2023}]\label{thm:comparison}
Let $M$ be a compact orientable Riemannian $n$-manifold with boundary components $\Gamma_1$, $\Gamma_0$. Suppose there exists a $\C^{1,1}$ function $u\colon M\to [0,1]$ with $\nabla u\neq 0$ on $M$, and $u=i$ on $\Gamma_i$. Let $\kappa:=(\kappa_1,\dots,\kappa_{n-1})$ be principal curvatures of level sets of $u$ with respect to $e_n:=\nabla u/|\nabla u|$, and let  $e_1,\dots, e_{n-1}$ be an orthonormal set of the corresponding principal directions. Then, 
for $0\leq r\leq n-1$,
\begin{multline}\label{eq:main}
\M_{r}(\G_1)-\M_{r}(\G_0)
=
(r+1)\int_{M}\sigma_{r+1}(\kappa)\, \\
+
\int_{M}\left(-\sum\kappa_{i_1}\dots\kappa_{i_{r-1}}K_{i_{r}n}
+
\frac{1}{|\nabla u|}\sum\kappa_{i_1}\dots\kappa_{i_{r-2}}|\nabla u|_{i_{r-1}}R_{i_{r}i_{r-1}i_{r}n}\right) ,
\end{multline}
where $|\nabla u|_i:=\nabla_{e_i} |\nabla u|$, $R_{ijkl}=\langle R(e_i, e_j)e_k, e_l\rangle$ are components of the Riemann curvature tensor of $M$, $K_{ij}=R_{ijij}$ is the sectional curvature, and
the sums range over distinct values of $1\leq i_1,\dots, i_r\leq n-1$,  
with $i_1<\dots <i_{r-1}$ in the first sum, and $i_1<\dots <i_{r-2}$ in the second sum. 
\end{theorem}

In \cite{ghomi-spruck2023}, the above theorem was established  via Reilly's identities \cite{reilly1977}. Here we present a  somewhat shorter and conceptually simpler proof using  differential forms  which we construct after Chern \cite{chern1945}, as Borb\'{e}ly \cite{borbely2002,borbely2002b} had also done earlier. More specifically, we devise a differential $(n-1)$-form  $\Phi_r$ on $M$ so that $\mathcal{M}_r(\Gamma_i)$ correspond to integration of $\Phi_r$ on $\Gamma_i$. Then computing the exterior derivative $d\Phi_r$ yields \eqref{eq:main} via Stokes theorem. Various applications of Theorem \ref{thm:comparison}  are developed in \cite{ghomi-spruck2022,  ghomi-spruck2023rigidity}, including  total curvature bounds,  and rigidity results in Riemannian geometry. See also \cite{ghomi-spruck2023} for more results of this type.

\section{Basic Formulas}
As in the statement of Theorem \ref{thm:comparison}, we let $M$ be a  compact orientable Riemannian $n$-manifold with boundary $\partial M=\Gamma_1\cup\Gamma_0$. Furthermore, $\langle\cdot,\cdot\rangle$ denotes the metric on $M$, with induced norm $|\cdot|:=\langle\cdot,\cdot\rangle^{1/2}$,  connection $\nabla$, and curvature operator
$$
R(X,Y)Z:=\nabla_Y\nabla_XZ-\nabla_X\nabla_YZ+\nabla_{[X,Y]}Z,
$$
for vector fields $X$, $Y$, $Z$ on $M$. The sectional curvature of $M$ with respect to a pair of orthonormal vectors $x$, $y$ in the tangent space $T_p M$ may be defined as
$$
K(x,y):=\langle R(X,Y)X, Y\rangle,
$$
where $X$, $Y$ are local extensions of $x$, $y$.
With $u$ as in the statement of Theorem \ref{thm:comparison}, and
for $0\leq t\leq 1$, let $\Gamma_t:=u^{-1}(t)$ be the level hypersurface of $u$ at height $t$. 
Since $u$ is $\C^{1,1}$, $\Gamma_t$ is twice differentiable almost everywhere by Rademacher's theorem.
At every such point $p$ of $\Gamma_t$,  let $e_i$, $i=1,\dots,n$, be the orthonormal frame mentioned above, i.e.,  
$$
e_n:=\frac{\nabla u}{|\nabla u|},
$$ 
and $e_1,\dots, e_{n-1}$ form a set of orthonormal principal directions of $\Gamma_t$ at $p$. Furthermore we assume that $e_i$ is \emph{positively oriented}, i.e., 
\begin{equation}\label{eq:dvol}
d\textup{vol}_M(e_1,\dots,e_n)=1,
\end{equation}
 where $d\textup{vol}_M$ denotes the volume form of $M$. We call $e_i$ a \emph{principal frame} associated to (level sets of) $u$.
Let $\theta^i$ be the corresponding dual one forms on $T_pM$ given by
\begin{equation}\label{eq:kronecker}
\theta^i(e_j)=\delta^i_j,
\end{equation}
where $\delta^i_j$ is the Kronecker  function.
 Note that $e_i$ may be extended to a $\C^1$ orthonormal frame $\ol e_i$ in a neighborhood of $p$ in $M$ so that $\ol e_n=e_n$  and thus $\ol e_1, \dots, \ol e_{n-1}$ remain tangent to $\Gamma_t$ (though they may no longer be principal directions).  The corresponding connection $1$-forms on $T_pM$ are then given by
\begin{equation*}\label{eq:skew-symmetry}
\omega^i_j(\cdot):=\langle\nabla_{(\cdot)} \ol e_j, e_i\rangle=-\langle  e_j, \nabla_{(\cdot)} \ol e_i\rangle=-\omega^j_i(\cdot),
\end{equation*}
for $1\leq i, j\leq n$.  Since $e_i$, $i=1,\dots, n-1$ are principal directions, and $\ol e_n=e_n$ is the normal of $\Gamma_t$, 
\begin{equation}\label{eq:principal}
\omega_n^i(e_j)=\langle\nabla_{e_j} \ol e_n, e_i\rangle=\delta^i_j\kappa_i, \quad\text{$1\leq i,j\leq n-1$},
\end{equation}
 where $\kappa_i$ are
the principal curvatures of $\Gamma_t$ with respect to $e_n$. We also record that, 
\begin{equation}\label{eq:en}
\omega_n^i(e_n)=\frac{\langle\nabla_{e_n} \nabla u,e_i\rangle}{|\nabla u|}=\frac{\langle\nabla_{e_i} \nabla u,e_n\rangle}{|\nabla u|}=\frac{\langle\nabla_{e_i} \nabla u,\nabla u\rangle}{|\nabla u|^2}=\frac{|\nabla u|_i}{|\nabla u|}, \quad\text{$1\leq i\leq n-1$},
\end{equation}
where $|\nabla u|_i=\nabla_{e_i}|\nabla u|$, and the second equality is due to the symmetry of the Hessian of $u$.
Next, we compute $\omega_i^j$ for $i,j\neq n$.
We may assume that $\ol e_1, \dots, \ol e_{n-1}$ are parallel translations of $e_1,\dots, e_{n-1}$ on $\Gamma_t$, i.e., $\ol\nabla_{e_i}\ol e_j=0$, for $1\leq i, j\leq n-1$ where $\ol\nabla:=\nabla^{\top}$ is the induced connection on $\Gamma_t$. Then $\omega_i^j(e_k)=\langle \ol\nabla_{e_k} \ol e_j, e_i\rangle=0$, for $1\leq i, j,k\leq n-1$.
Furthermore, we may assume that $\ol e_1, \dots, \ol e_{n-1}$ are parallel translated along the integral curve of $e_n$. Then $\nabla_{e_n}\ol e_i=0$ for $1\leq i\leq n-1$, which yields $\omega_i^j(e_n)=0$, for $1\leq i,j\leq n-1$. So we record that
\begin{equation}\label{eq:omegaij}
\omega_i^j=0, \quad 1\leq i,j\leq n-1.
\end{equation}
Cartan's structure equations state that
\begin{equation}\label{eq:cartan}
d\theta^i=\sum_{j=1}^n\theta^j\wedge\omega^i_j
 \quad\quad\;\;\text{and}\quad \quad\;\;
d\omega_j^i=\Omega_j^i-\sum_{k=1}^n\omega^k_j\wedge\omega^i_k,
\end{equation}
where $\Omega_j^i$ are the curvature $2$-forms given by
$$
\Omega_j^i(e_\ell, e_k):=-\big\langle R(e_\ell, e_k)e_j, e_i\big\rangle=\big\langle R(e_\ell, e_k)e_i, e_j\big\rangle=:R_{\ell k i j}.
$$
Note that $R_{\ell k i j}=-R_{k\ell i j}$. We also set 
\begin{equation}\label{eq:K}
K_{ij}:=K(e_i, e_j)=R_{ijij}.
\end{equation}
 Finally we record some basic formulas from exterior algebra which will be used in the next section. If $\lambda$ is a $k$-form,  and $\phi$ is an $\ell$-form, then 
\begin{equation}\label{eq:theta-phi}
\lambda\wedge\phi(e_1,\dots, e_{k+\ell})=\sum\epsilon(i_1\dots i_{k+\ell})\,\lambda(e_{i_1},\dots, e_{i_k})\,\phi(e_{i_{k+1}},\dots,e_{i_{k+\ell}})
\end{equation}
where the sum ranges over  $1\leq i_1,\dots ,i_{k+\ell}\leq k+\ell$, with $i_1<\dots< i_k$, and $i_{k+1}<\dots <i_{k+\ell}$; furthermore, $\epsilon(i_1\dots i_n):=1$, or $-1$ depending on whether $i_1\dots i_n$ is an even or odd permutation of $1\dots n$ respectively. Note that
\begin{equation}\label{eq:epsilon}
\epsilon(i_1 \dots i_{r-1}n i_{r+1}\dots i_{n-1})=(-1)^{n-1-r} \epsilon(i_1 \dots i_{n-1}),
\end{equation}
since $\epsilon(i_1 \dots i_{n-1})=\epsilon(i_1 \dots i_{n-1}n)$.
The following  identities will also be useful
\begin{gather}\label{eq:d}
d(\theta^1\wedge\dots\wedge\theta^k)
=
\sum \epsilon(i_1\dots i_k)\,d\theta^{i_1}\wedge\theta^{i_2}\wedge\dots\wedge \theta^{i_k}\\ \notag
=
(-1)^{k-1}\sum \epsilon(i_1\dots i_k)\,\theta^{i_1}\wedge\dots \wedge\theta^{i_{k-1}}\wedge d\theta^{i_k},
\end{gather}
where the sums range over  $1\leq i_1,\dots, i_k\leq k$ with $i_2<\dots<i_{k}$ in the first sum, and $i_1<\dots<i_{k-1}$ in the second sum. 

\section{Proof of Theorem \ref{thm:comparison}}
Let $\theta^i$ be the dual $1$-forms, and $\omega^i_j$ be the connection forms corresponding to the principal frame $e_i$ of $u$ discussed in the last section. For $0\leq r\leq n-1$, we define the $(n-1)$-forms
$$
\Phi_r:=\sum\epsilon(i_1\dots i_{n-1})\,\omega^{i_1}_{n}\wedge\dots\wedge\omega^{i_r}_{n}\wedge\theta^{i_{r+1}}\wedge\dots\wedge\theta^{i_{n-1}},
$$
where the sum ranges over  $1\leq i_1,\dots,i_{n-1}\leq n-1$ with
$i_1<\dots<i_{r}$, and $i_{r+1}<\dots<i_{n-1}$.
For $r=n-1$, this form appears  in Chern \cite{chern1945}, and later in Borb\'{e}ly \cite{borbely2002} (where it is denoted 
as  ``$\Phi_0$'' and ``$\Phi$'' respectively). The form $\Phi_1$ has also been used by Borb\'{e}ly in \cite{borbely2002b}.
One quickly checks, using  \eqref{eq:kronecker}, \eqref{eq:principal}, and \eqref{eq:theta-phi}, that
\begin{eqnarray}\label{eq:phir}
\Phi_r(e_1,\dots, e_{n-1})=
\sigma_r(\kappa),
\end{eqnarray}
which is the main feature of these forms. Recall that $\Gamma_t:=u^{-1}(t)$ is the level hypersurface of $u$ at height $t$, for $0\leq t\leq 1$. Let $\Phi_r|_{\Gamma_t}$ denote the pull back of $\Phi_r$ via the inclusion map $\Gamma_t\to M$.
Since $\Phi_r|_{\Gamma_t}$ is an $(n-1)$-form on $\Gamma_t$, it is a multiple of the volume form of $\Gamma_t$, which is given by
\begin{equation}\label{eq:dvol-t}
d\textup{vol}_{\Gamma_t}(e_1,\dots, e_{n-1}):=d\textup{vol}_M(e_n,e_1,\dots, e_{n-1})=\epsilon(n 1\dots n-1)=(-1)^{n-1}.
\end{equation}
Note that here we have used the assumption \eqref{eq:dvol} that $e_i$ is positively oriented. So it follows from \eqref{eq:phir} and \eqref{eq:dvol-t} that 
\begin{equation}\label{eq:phi-Gt}
\Phi_{r}|_{\Gamma_t}=(-1)^{n-1}\sigma_{r}(\kappa) \,d\textup{vol}_{\Gamma_t}.
\end{equation}
This shows that $\Phi_r$ depends only on $e_n$, not the choice of $e_1,\dots, e_{n-1}$ (which also follows from transformation rules for $\omega^i_n$ and $\theta^i$ under a change of frame $e_i\to e_i'$ with $e_n=e_n'$; see \cite[p. 269]{borbely2002b}).
In addition, \eqref{eq:phi-Gt} shows that
$$
\M_r(\Gamma_t)=\int_{\Gamma_t}\sigma_r(\kappa):=\int_{\Gamma_t}\sigma_r(\kappa)\,d\textup{vol}_{\Gamma_t}=(-1)^{n-1}\int_{\Gamma_t}\Phi_r.
$$
Consequently, by Stokes theorem, for the left hand side of \eqref{eq:main} we have
\begin{equation}\label{eq:MGg}
\M_{r}(\G_1)-\M_{r}(\Gamma_0)
=
(-1)^{n-1}\int_{\partial M}\Phi_{r}
=
(-1)^{n-1}\int_{M}d\Phi_{r}.
\end{equation}
Here we have  used the assumption that $u|_{\Gamma_1}>u|_{\Gamma_0}$, which ensures that $e_n$ points outward on $\Gamma_1$ and inward on $\Gamma_0$ with respect to $M$. Furthermore, since $\Phi_r$ depends only on $e_n$ and $u$ is $\C^{1,1}$, it follows that $\Phi_r$ is Lipschitz (in local coordinates). Hence $d\Phi_r$ is integrable, and the use of Stokes theorem here is justified.

Next we compute $d\Phi_r$. Since $\omega^{i_1}_{n}\wedge\dots\wedge\omega^{i_r}_{n}$ is an $r$-form,  the product rule for exterior differentiation yields that
\begin{multline}\label{eq:dphir-0}
d\Phi_r
=
(-1)^{r} \sum\epsilon(i_1\dots i_{n-1})\,\omega^{i_1}_{n}\wedge\dots\wedge\omega^{i_r}_{n}\wedge d\big(\theta^{i_{r+1}}\wedge\dots\wedge\theta^{i_{n-1}}\big)\\
+
\sum \epsilon(i_1\dots i_{n-1})\,d\big(\omega^{i_1}_{n}\wedge\dots\wedge\omega^{i_r}_{n}\big)\wedge\theta^{i_{r+1}}\wedge\dots\wedge\theta^{i_{n-1}},
\end{multline}
where  the sums still range over $i_1<\dots <i_r$ and $i_{r+1}<\dots<i_{n-1}$.
By \eqref{eq:d},  the structure equations \eqref{eq:cartan}, and \eqref{eq:omegaij}, the first term in \eqref{eq:dphir-0} reduces to
\begin{gather*}
\;\;\;\;\,(-1)^{r+1}\sum \epsilon(i_{1}\dots i_{n-1})\, \omega^{i_1}_{n}\wedge\dots\wedge\omega^{i_r}_{n}\wedge 
\omega^{i_{r+1}}_n\wedge\theta^n\wedge\theta^{i_{r+2}}\wedge\dots\wedge\theta^{i_{n-1}}\\
= 
(-1)^{n-1}\sum \epsilon(i_1\dots i_{n-1}) \,\omega^{i_1}_{n}\wedge\dots\wedge\omega^{i_r}_{n}\wedge 
\omega^{i_{r+1}}_n\wedge\theta^{i_{r+2}}\wedge\dots\wedge\theta^{i_{n-1}}\wedge\theta^n\\
= 
(-1)^{n-1}(r+1)\Phi_{r+1}\wedge\theta^n,
\end{gather*}
where the sums now range over $i_1<\dots <i_r$, and $i_{r+2}<\dots <i_{n-1}$. The factor $(r+1)$ appears in the last line because definition of $\Phi_{r+1}$ requires that $i_1<\dots <i_{r+1}$. Applying \eqref{eq:d} and  \eqref{eq:cartan} also to the second term in \eqref{eq:dphir-0}, we obtain
\begin{multline}\label{eq:dphir}
d\Phi_r
=(-1)^{n-1}(r+1)\Phi_{r+1}\wedge\theta^{n}\\
+(-1)^{r-1}
\sum\epsilon(i_1\dots i_{n-1})\,\omega^{i_1}_n\wedge\dots\wedge\omega^{i_{r-1}}_{n}\wedge\Omega^{i_r}_{n}\wedge\theta^{i_{r+1}}\wedge\dots\wedge\theta^{i_{n-1}},
\end{multline}
where the sum ranges over $i_1<\dots <i_{r-1}$, and $i_{r+1}<\dots <i_{n-1}$.
For $r=1$, this formula had been computed earlier by Borb\'{e}ly  \cite[(6)]{borbely2002b}.

By \eqref{eq:MGg}, it remains to show that $(-1)^{n-1}\int_M d\Phi_r$ yields the right hand side of \eqref{eq:main}. To see this first note that, by \eqref{eq:theta-phi} and \eqref{eq:phir},
$$
\Phi_{r+1}\wedge\theta^{n}=\Phi_{r+1}\wedge\theta^{n}(e_1,\dots,e_n)\,d\textup{vol}_M=\sigma_{r+1}(\kappa)\,d\textup{vol}_M.
$$
Thus the first term on the right hand side of  \eqref{eq:dphir} quickly yields the first integral on the right hand side of \eqref{eq:main}. To obtain the second  integral there, we evaluate the sum in  \eqref{eq:dphir} at $e_i$, which yields
\begin{gather*}\label{eq:lastsum}
\sum\epsilon(j_1\dots j_n)\,\epsilon(i_1\dots i_{n-1})\,\omega^{i_1}_n(e_{j_1})\dots\omega^{i_{r-1}}_{n}(e_{j_{r-1}})\Omega^{i_r}_{n}(e_{j_r}, e_{j_{r+1}})\theta^{i_{r+1}}(e_{j_{r+2}})\dots\theta^{i_{n-1}}(e_{j_n})\\ \notag
=
\sum\epsilon(j_1\dots j_{r+1}i_{r+1}\dots i_{n-1})\,\epsilon(i_1\dots i_{n-1})\,\omega^{i_1}_n(e_{j_1})\dots\omega^{i_{r-1}}_{n}(e_{j_{r-1}})R_{j_r j_{r+1}i_rn},
\end{gather*}
where  the sums range over $1\leq j_1\dots j_n\leq n$ with $j_r<j_{r+1}$ by \eqref{eq:theta-phi}, and the range for $1\leq i_1,\dots, i_{n-1}\leq n-1$ remains as in \eqref{eq:dphir}, i.e., $i_1<\dots <i_{r-1}$, and $i_{r+1}<\dots <i_{n-1}$.
The last sum may be partitioned into $A+B$, where $A$ consists of terms with $j_{r+1}=n$, and $B$ of terms with $j_{r+1}\neq n$. If $j_{r+1}=n$, then $j_1,\dots, j_{r-1}\neq n$, which yields $i_k=j_{k+1}$ for $k=1,\dots, r-2$ by \eqref{eq:principal}. This in turn forces $j_r=i_r$, as they are the only remaining indices. So by \eqref{eq:epsilon} and \eqref{eq:K}, 
\begin{gather*}
A
=
\sum\epsilon(i_1 \dots i_{r-1}n i_{r+1}\dots i_{n-1})\,\epsilon(i_1\dots i_{n-1})\,\kappa_{i_1}\dots\kappa_{i_{r-1}}R_{i_r n  i_rn}\\
=
(-1)^{n-r-1}\sum\kappa_{i_1}\dots\kappa_{i_{r-1}}K_{i_rn},
\end{gather*}
where we still have $i_1<\dots<i_{r-1}$.
This yields the first term in the second integral in \eqref{eq:main}, after multiplication by the sign factors $(-1)^{r-1}$ from \eqref{eq:dphir} and $(-1)^{n-1}$ from \eqref{eq:MGg}, which ensures the desired sign  $-1$.
Next, to compute $B$, note that if $j_{r+1}\neq n$, then $j_r\neq n$ either, since $j_r<j_{r+1}$, which forces $j_k=n$, for some $1\leq k\leq r-1$.
We may assume $k=r-1$ after reindexing. Then $j_1,\dots, j_{r-2}\neq n$, which yields $i_k=j_{k}$ for $k=1,\dots, r-2$ by \eqref{eq:principal}.
So by \eqref{eq:en}
\begin{gather*}
B
=
\sum\epsilon(i_1\dots  i_{r-2}n j_r j_{r+1}i_{r+1}\dots i_{n-1})\,\epsilon(i_1\dots i_{n-1})\,\kappa_{i_1}\dots\kappa_{i_{r-2}}\frac{|\nabla u|_{i_{r-1}}}{|\nabla u|}R_{j_rj_{r+1} i_rn}\\
=
\sum\epsilon(i_1\dots  i_{r-2}n i_{r-1}\dots i_{n-1})\,\epsilon(i_1\dots i_{n-1})\,\kappa_{i_1}\dots\kappa_{i_{r-2}}\frac{|\nabla u|_{i_{r-1}}}{|\nabla u|}R_{i_{r-1}i_{r} i_rn}\\
=
(-1)^{n-r}\sum\kappa_{i_1}\dots\kappa_{i_{r-2}}\frac{|\nabla u|_{i_{r-1}}}{|\nabla u|}R_{i_{r}i_{r-1} i_r n},
\end{gather*}
where the second equality  holds because $\{j_r,j_{r+1}\}=\{i_{r-1},i_r\}$, since these are the only remaining indices. We may assume then that
$j_r=i_{r-1}$, and $j_{r+1}=i_r$, since switching $j_r$ and $j_{r+1}$ does not change the sign of the right hand side of the first equality for $B$.
The sign $(-1)^{n-r}$ in the third equality is due to \eqref{eq:epsilon} and  switching two  indices in the Riemann tensor coefficient. 
Finally note that the restriction on the range of indices in the last sum is now $i_1<\dots <i_{r-2}$, since $i_{r-1}$ corresponds to $j_{r-1}$, and we set $r-1=k$ during the reindexing above.
So $B$ yields  the second term in the second integral in \eqref{eq:main}, after multiplication by  $(-1)^{r-1}$ and $(-1)^{n-1}$, as was the case for $A$, which ensures the desired sign $+1$. This concludes the proof of Theorem \ref{thm:comparison}.

\addtocontents{toc}{\protect\setcounter{tocdepth}{0}}

\section*{Acknowledgment}
This work is an outgrowth of  extensive collaborations with Joel Spruck on the topic of total curvature, and is indebted to him for numerous discussions.

\addtocontents{toc}{\protect\setcounter{tocdepth}{1}}

\bibliography{references}

\begin{bibdiv}
\begin{biblist}

\bib{borbely2002b}{article}{
      author={Borb\'{e}ly, Albert},
       title={An estimate for the total mean curvature in negatively curved
  spaces},
        date={2002},
        ISSN={0004-9727},
     journal={Bull. Austral. Math. Soc.},
      volume={66},
      number={2},
       pages={267\ndash 273},
         url={https://doi.org/10.1017/S0004972700040119},
      review={\MR{1932350}},
}

\bib{borbely2002}{article}{
      author={Borb\'{e}ly, Albert},
       title={On the total curvature of convex hypersurfaces in hyperbolic
  spaces},
        date={2002},
        ISSN={0002-9939},
     journal={Proc. Amer. Math. Soc.},
      volume={130},
      number={3},
       pages={849\ndash 854},
  url={https://doi-org.prx.library.gatech.edu/10.1090/S0002-9939-01-06101-9},
      review={\MR{1866041}},
}

\bib{chern1945}{article}{
      author={Chern, Shiing-shen},
       title={On the curvatura integra in a {R}iemannian manifold},
        date={1945},
        ISSN={0003-486X},
     journal={Ann. of Math. (2)},
      volume={46},
       pages={674\ndash 684},
         url={https://doi-org.prx.library.gatech.edu/10.2307/1969203},
      review={\MR{14760}},
}

\bib{ghomi-spruck2022}{article}{
      author={Ghomi, Mohammad},
      author={Spruck, Joel},
       title={Total curvature and the isoperimetric inequality in
  {C}artan-{H}adamard manifolds},
        date={2022},
        ISSN={1050-6926},
     journal={J. Geom. Anal.},
      volume={32},
      number={2},
       pages={Paper No. 50, 54},
         url={https://doi.org/10.1007/s12220-021-00801-2},
      review={\MR{4358702}},
}

\bib{ghomi-spruck2023rigidity}{article}{
      author={Ghomi, Mohammad},
      author={Spruck, Joel},
       title={Rigidity of nonpositively curved manifolds with convex boundary},
        date={2023},
     journal={arXiv:2210.05588, to appear in Proc. Amer. Math. Soc.},
}

\bib{ghomi-spruck2023}{article}{
      author={Ghomi, Mohammad},
      author={Spruck, Joel},
       title={Total mean curvatures of {R}iemannian hypersurfaces},
        date={2023},
     journal={Advanced Nonlinear Studies},
      volume={23},
      number={1},
       pages={20220029},
         url={https://doi.org/10.1515/ans-2022-0029},
}

\bib{reilly1977}{article}{
      author={Reilly, Robert~C.},
       title={Applications of the {H}essian operator in a {R}iemannian
  manifold},
        date={1977},
        ISSN={0022-2518},
     journal={Indiana Univ. Math. J.},
      volume={26},
      number={3},
       pages={459\ndash 472},
  url={https://doi-org.prx.library.gatech.edu/10.1512/iumj.1977.26.26036},
      review={\MR{0474149}},
}

\end{biblist}
\end{bibdiv}

\end{document}